# An Upper Bound for Lebesgue's Covering Problem

Philip Gibbs

Abstract: A covering problem posed by Henri Lebesgue in 1914 seeks to find the convex shape of smallest area that contains a subset congruent to any point set of unit diameter in the Euclidean plane. Methods used previously to construct such a covering can be refined and extended to provide an improved upper bound for the optimal area. An upper bound of 0.8440935944 is found.

## Introduction and history

Henri Lebesgue is known for his theory of integration based on measure theory. He defined the measure of a point set as the minimum area of covering sets whose area can be calculated in a more conventional way. In 1914 he wrote a letter to his colleague Julius Pál in which he asked about a related problem [1]. What is the smallest area $a$ of a convex planar shape that can cover any point set of unit diameter in the 2D plane? It is understood that the set can be translated, rotated and reflected to fall inside the covering. The diameter of a bounded pint set is defined as the supremum of the distances between pairs of points in the set.

Pál set about trying to solve the problem. A prior result by Heinrich Jung had shown that any planar shape of diameter one can be covered by a circle of radius $\frac{1}{\sqrt{3}}$. This sets a weak upper bound for the optimal area $a \leq \frac{\pi}{3}$ [2,3]. In fact it is easy to improve on this because a shape of diameter 1 is easily seen to fit in a square of unit side using translations alone, so $a \leq 1$. A simple lower bound $a \geq \frac{\pi}{4}$ comes from the need to cover a circle of diameter one. Pál studied the problem further and published a paper in 1920 with improved upper and lower bounds [4]. First he demonstrated that a shape of unit diameter can be covered by a regular hexagon with unit distance between opposite sides, i.e. the regular hexagon inscribed inside Jung's circular covering. He then showed that two of the corners can be removed from this hexagon along the edge of the largest regular dodecahedron that fits inside the hexagon. This reduced the upper bound to $a \leq 2 - \frac{2}{\sqrt{3}} = 0.84529946 \dots$ He also found the area of the smallest shape covering a circle and equilateral triangle to set an improved lower bound of $a \geq \frac{\pi}{8} + \frac{3}{4} = 0.82571178 \dots$

In 1936 Roland Sprague pointed out that an area of Pál's covering near one corner bounded by two arcs could be further removed [5]. The area of the remaining covering was later calculated to be $0.84413770 \dots$ [6] This was thought to be the best possible answer until Hansen shows that some imperceptibly small regions of total area about $4 \times 10^{-11}$ could still be removed [7,8,9]. Hansen's reduction was the first that used reflection rather than just translation and rotation. Setting better lower bounds has been harder. György Elekes improved the lower bound to 0.8271 by combining polygons with sides a power of three [10]. Then Brass and Sharifi raised it to 0.832 by combining the circle and triangle with a regular pentagon [11]. Further progress on the upper bound looked difficult until a computational study revealed that an improvement could be made by cutting off Pál's corners at a different angle [12]. In 2015, this led Baez, Bagdasaryan and Gibbs to prove a new upper bound of 0.8441153 [13]. It was known at that time that some further small reductions are possible at the cost of increased complexity. In the present work the same methods are taken to their logical conclusion to give a new upper bound.

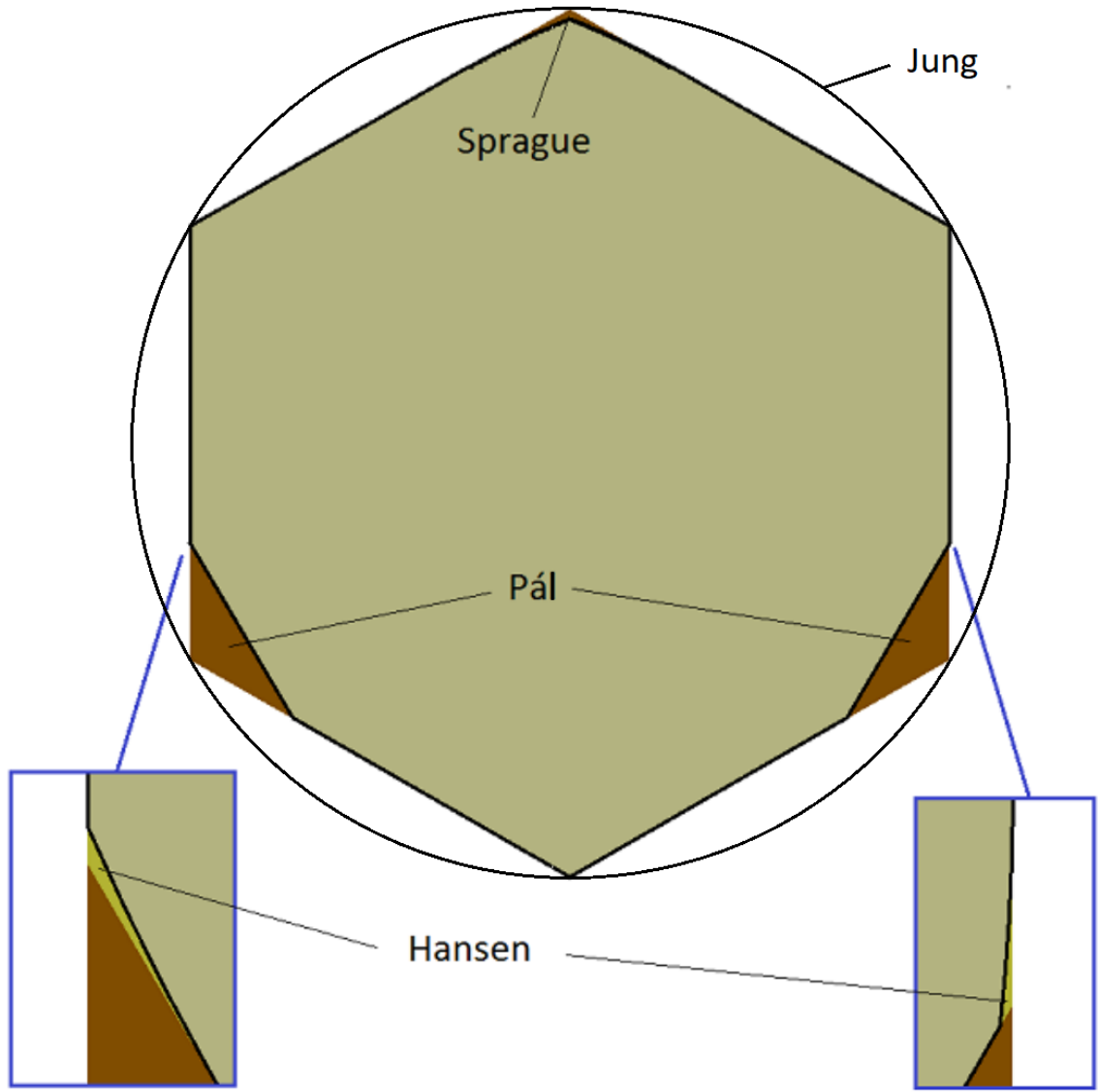

**Figure 1**

## Preliminaries and definitions

The **diameter** of a bounded point set in a metric space is by definition the supremum of all distances between pairs of points in the set. Two point sets in a metric space are **isometric** or **congruent** if there is a one-to-one mapping between their points such that the distance between any two points is equal to the distance between the mappings of the two points.

**Definition 1**: A **covering** for a collection of point sets in a metric space is a set that contains a subset congruent to each point set in the collection.

In the plane a covering for a collection of shapes is a shape that can cover a copy of every shape in the collection when it is transformed by some combination of translations, rotations and reflections.

**Lebesgue's covering problem** seeks the minimum area for a convex covering of all point sets of unit diameter in the Euclidean plane. Sometimes the convexity requirement is dropped. It can be shown that a shape (not necessarily unique) exists with the minimum area. This was demonstrated for the convex case by Kelly and Weiss in 1979 using the Blaschke selection theorem [14] and for the non-convex case by Kovalev in 1985, who also demonstrated that the shape is a star-domain [15]. Here only the convex case is considered further.

If a convex covering is an open set in the topological sense, it can be replaced by its closure without increasing its area. There is therefore no loss of generality in considering only closed shapes as covers and this will be assumed throughout. In particular, where a region is removed from the covering it is implicit that the boundary is left so that the covering remains closed.

**Definition 2**: A **maximal set of diameter $R$** in a metric space is one which is not the proper subset of any other set of diameter $R$. An **orbiform of diameter $R$** is a maximal set of diameter $R$ in the Euclidean 2D plane.

The boundary of an orbiform is called a **curve of constant width** because its width measured between any pair of parallel lines touching the shape on either side is always the same as its diameter. Any bounded point set in a metric space is the subset of some maximal set of the same diameter. This follows from Zorn's lemma with set inclusion as the partial ordering. In the case of the 2D Euclidean plane we use the following proposition.

**Proposition 1**: Any bounded point set of unit diameter in the plane is a subset of an orbiform of unit diameter. Proof: see Vrecica (1981) [16].

**Corollary 1**: A shape in the plane is a covering for all point sets of unit diameter if and only if it is a covering for all orbiforms of unit diameter.

A circle is the simplest example of a curve of constant width, with the disk as the corresponding orbiform, but other shapes with the same property are possible. While the circle has the largest enclosed area for a given diameter, the Rouleaux triangle has the least. Between these extremes are many other irregular shapes.

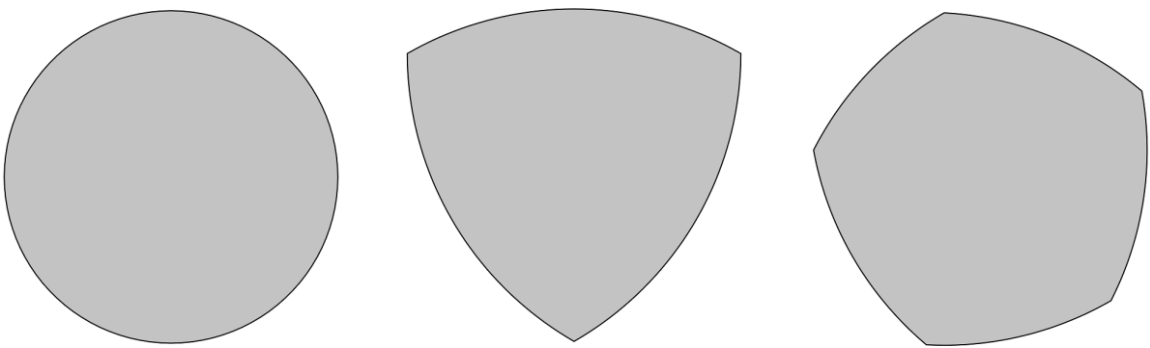

**Figure 2**

## Slanted Pál coverings

In 1920 Pál demonstrated that a hexagon with two corners removed is a covering [4]. In 2015 Baez, Bagdasaryan and Gibbs looked at a modified version of this covering [13].

**Definition 3**: The **slanted Pál covering** $\mathcal{P}(\sigma)$ for a slant angle $0 \leq \sigma < 30°$ is defined as follows

Start from a regular hexagon with unit distance between opposite sides. Label the corners clockwise from the top $A, B, C, D, E, F$. A copy of the hexagon is rotated through an angle of $30° + \sigma$ about its centre, then two of its sides are used to cut off the corners of the original hexagon at $C$ and $E$. The remaining interior and boundary is the shape $\mathcal{P}(\sigma)$ dependent on the slant angle $0 \leq \sigma < 30°$ (figure 3).

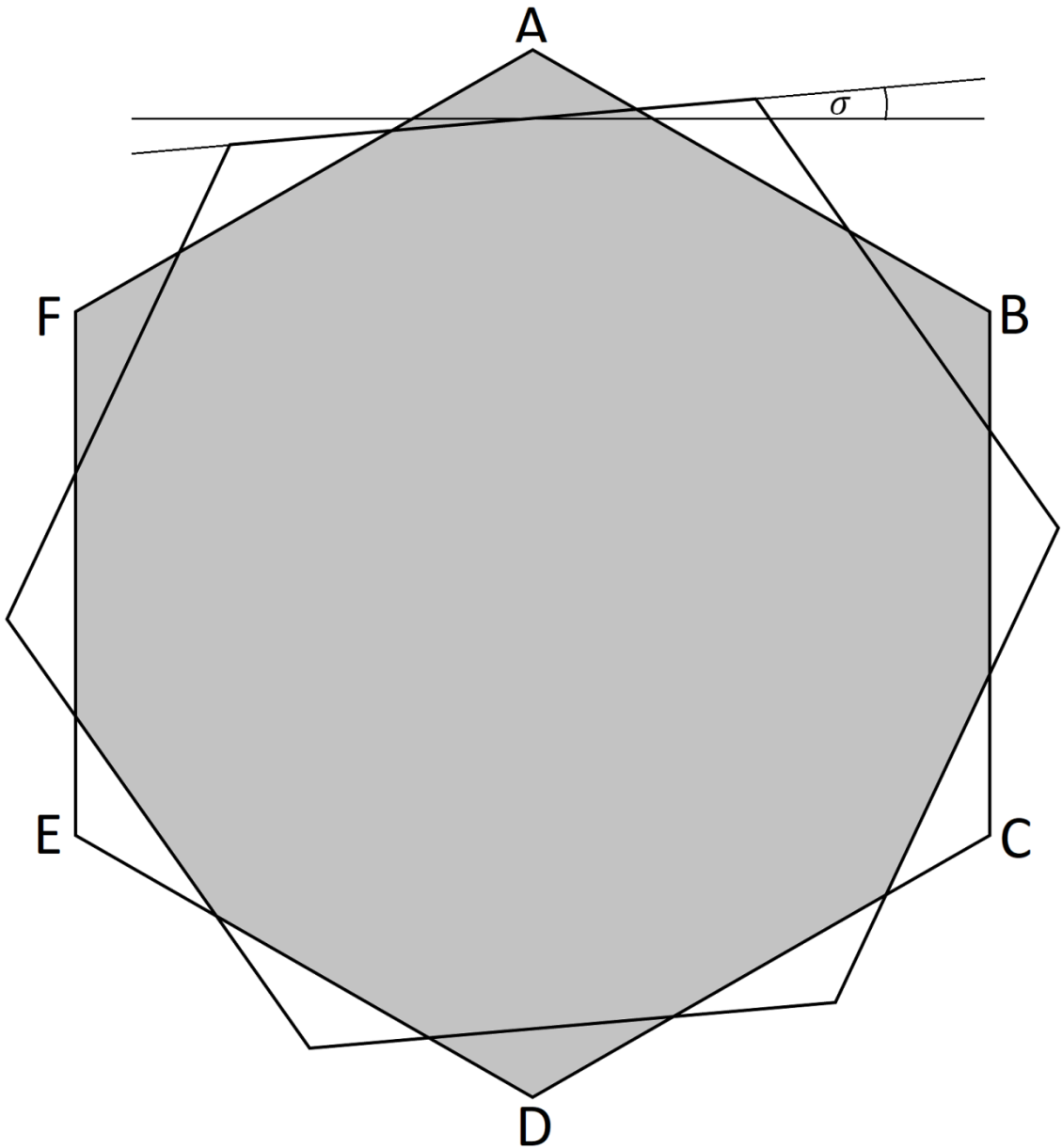

**Figure 3**

**Proposition 2**: $\mathcal{P}(\sigma)$ for a slant angle $0 \leq \sigma < 30°$ is a covering for all sets of unit diameter. Proof: See Baez, Bagdasaryan and Gibbs (2015) [13]

The minimum area of $\mathcal{P}(\sigma)$ is Pál's covering at $\sigma = 0$. However, there are further regions that can be removed to provide a smaller covering for non-zero slant angles. Computation suggests that the covering of minimum area is indeed a subset of $\mathcal{P}(\sigma)$ for a small value of $\sigma$ between one and two degrees [12]. A complete solution to Lebesgue's covering problem would require a proof of this which is currently unobtainable. In this work coverings that sit inside $\mathcal{P}(\sigma)$ are sought.

## Sprague-like reductions

Sprague was able to show that an area of Pál's original covering $\mathcal{P}(0)$ can be removed because any curve of constant width fitted into the covering cannot enter the region. Similar arguments can be applied to $\mathcal{P}(\sigma)$. In fact, for a non-zero slant angle more regions can be removed by similar arguments.

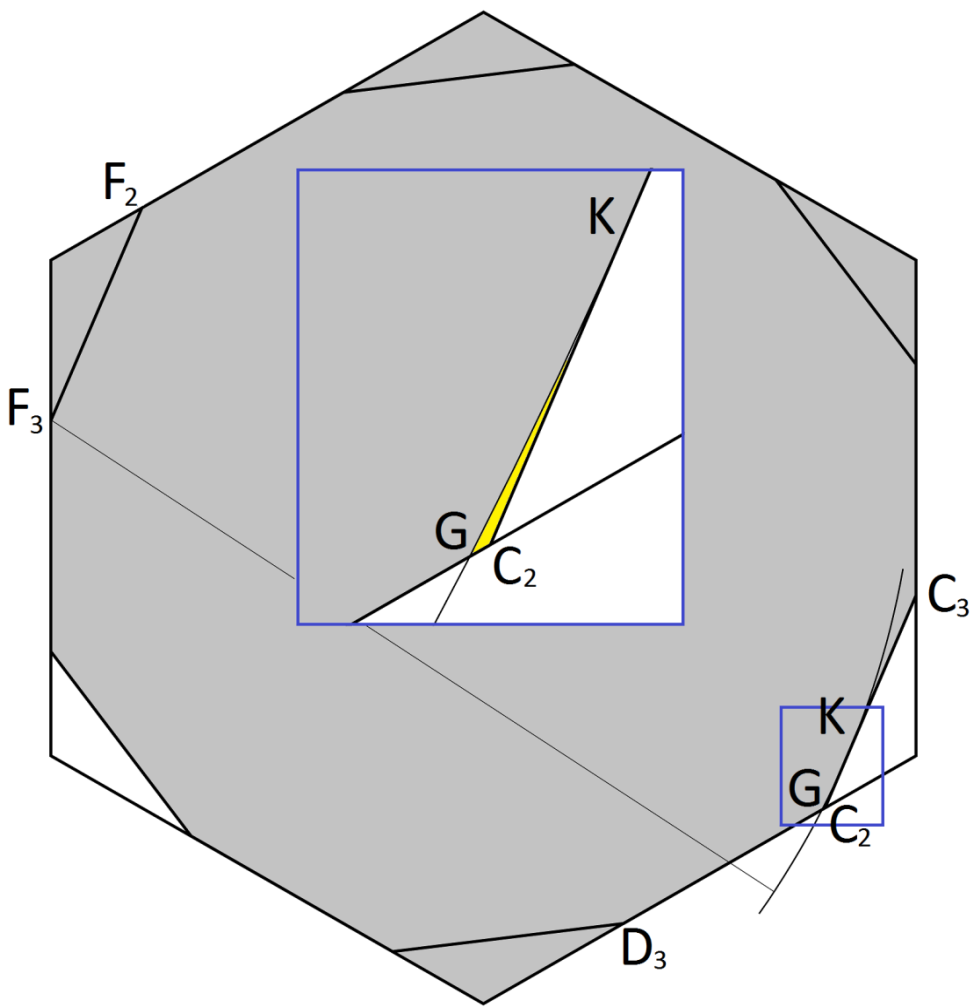

**Figure 4**

An orbiform of unit diameter that is a subset of $\mathcal{P}(\sigma)$ will not enter the interior of the removed corner at $C$. This means it must touch or cross the line segment $F_3F_2$ (see Figure 4.) All points in the orbiform must therefore be within unit distance of some point on this line segment. An arc of radius one centred on $F_3$ will touch the line segment $C_3C_2$ at a point $K$ and will cut the line segment $D_3C_2$ at a point $G$.

**Definition 4**: The **region $C_S$** is the closed shape $KGC_2$ bounded by the arc $KG$ and the two straight line segments $GC_2$ and $C_2K$

**Proposition 3**: If an orbiform $\mathcal{O}$ of unit diameter is a subset of $\mathcal{P}(\sigma)$ then no point of $\mathcal{O}$ is in the interior of the region $C_S$.

Proof: $\mathcal{O}$ must contain a point on the line segment $F_3F_2$ but all points in the interior of the region $C_S$ will be further than unit distance from that point and therefore cannot be in $\mathcal{O}$.

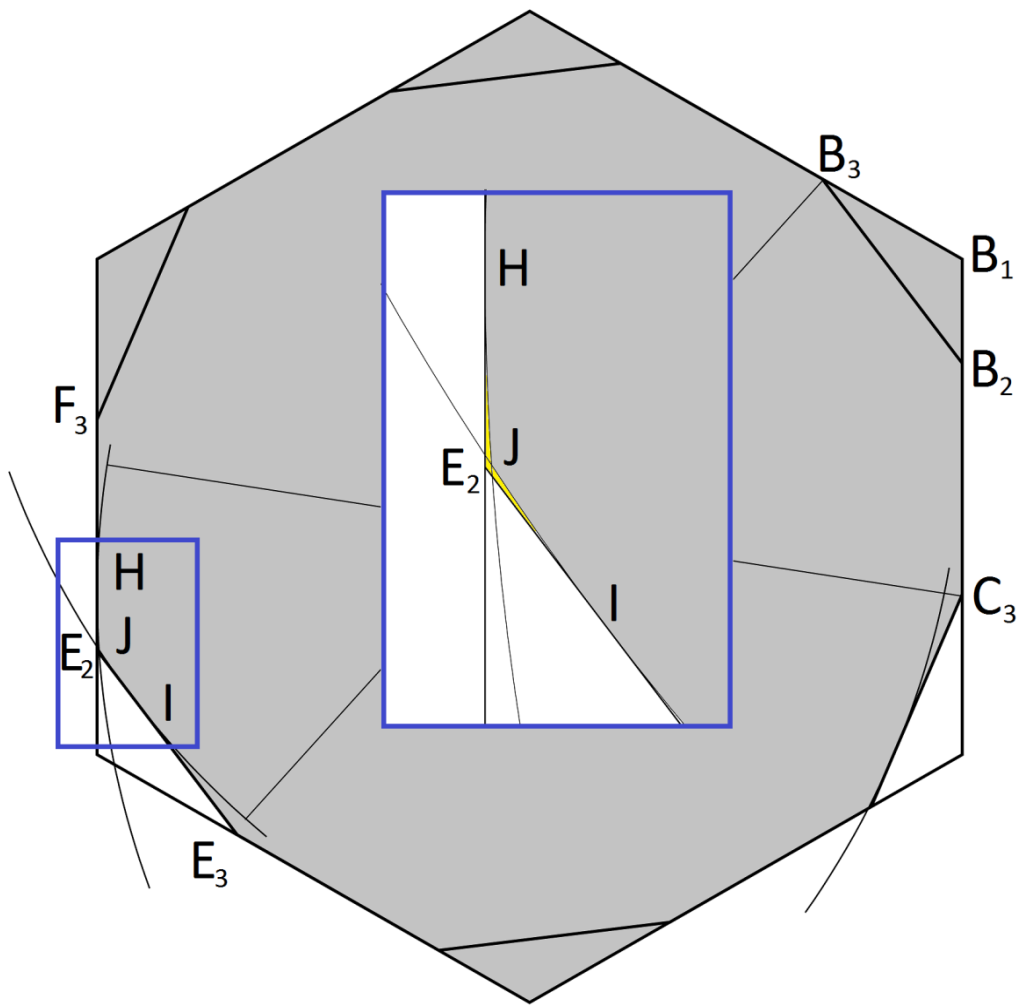

**Figure 5**

An orbiform of unit diameter fitted in $\mathcal{P}(\sigma)$ must also touch a point on the line segment $B_1C_3$ and it must touch or cross the line segment $B_3B_2$. An arc of unit radius centred on $B_3$ touches the line segment $E_3E_2$ at a point $I$. An arc of unit radius centred on $C_3$ touches the line segment $E_2F_3$ at a point $H$. These two arc segments intersect at a point $J$ (see Figure 5.)

**Definition 5**: The **region $E_S$** is the closed shape $IJHE_2$ bounded by the two arcs and two straight line segments.

**Proposition 4**: If an orbiform $\mathcal{O}$ of unit diameter is a subset of $\mathcal{P}(\sigma)$ then no point of $\mathcal{O}$ is in the interior of the region $E_S$.

Proof: All points in the interior of $E_S$ are at greater than unit distance from all points on $B_1C_3$ or all points on $B_3B_2$. $\mathcal{O}$ must contain points on both these line segments. Therefore no point in the interior of $E_S$ can be in $\mathcal{O}$.

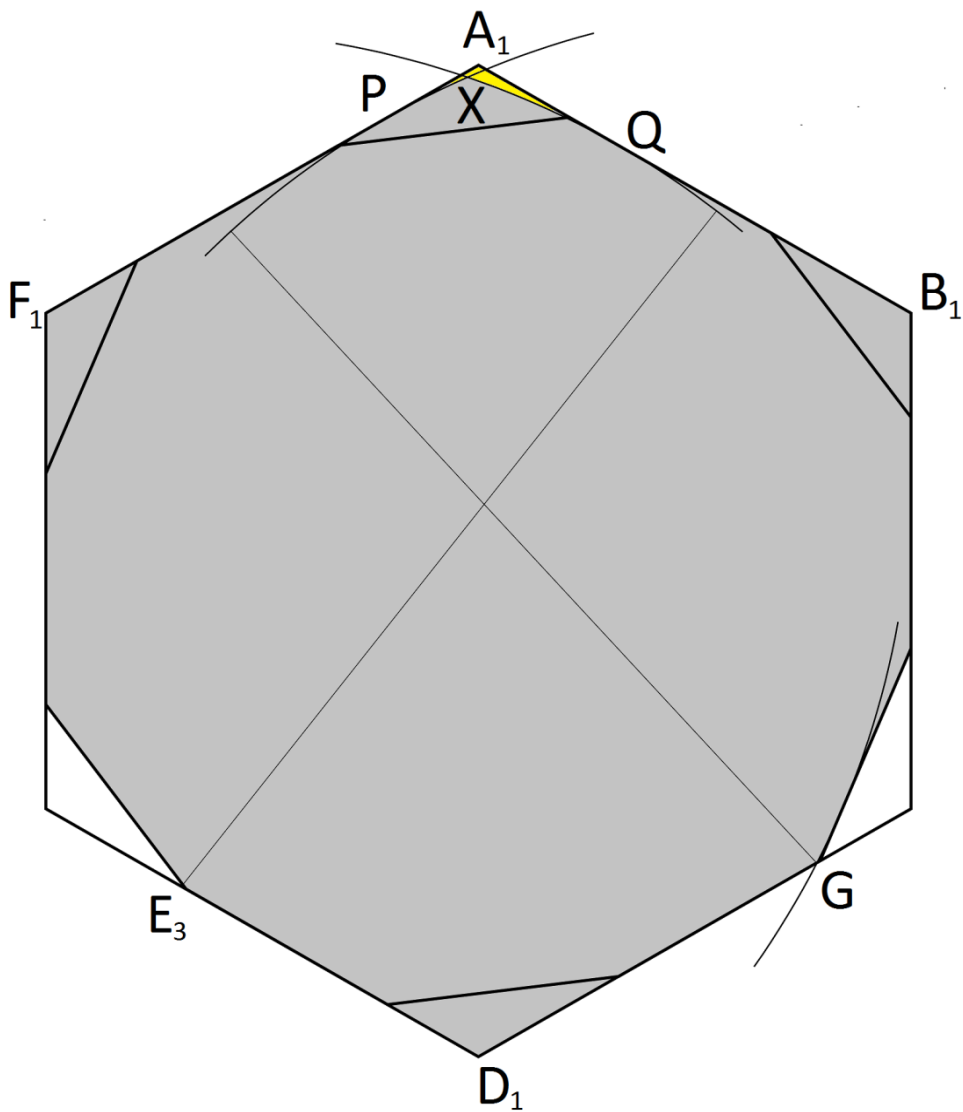

**Figure 6**

One further Sprague-like argument can be used. An orbiform of unit diamter fitted into $\mathcal{P}(\sigma)$ must touch the line segments $GD_1$ and $D_1E_3$. An arc of radius one centred on $G$ will touch $F_1A_1$ at a point $P$ and an arc of radius one centred on $E_3$ will touch $A_1B_1$ at a point $Q$. The two arcs will intersect at a point $X$. (See figure 6.)

**Definition 6**: The region $A_S$ is the closed shape $QXPA_1$ bounded by the two arcs and two straight line segments

**Proposition 5**: If an orbiform $\mathcal{O}$ of unit diameter is a subset of $\mathcal{P}(\sigma)$ then no point of $\mathcal{O}$ is in the interior of the region $A_S$.

Proof: All points in the interior of $A_S$ are at a distance greater than one from all points on $D_1E_3$ or at a distance greater than one from all points on $GD_1$. However, $\mathcal{O}$ must include points on both these line segments.

**Definition 7**: The **shape** $\mathcal{S}(\sigma)$ for $0 \leq \sigma < 30°$ is formed from the shape $\mathcal{P}(\sigma)$ (definition 3) by removing the regions $C_S$, $E_S$ and $A_S$ (definition 4,5 and 6) and forming the closure.

**Proposition 6**: The shape $\mathcal{S}(\sigma)$ is a covering for shapes of unit diameter.

**Lemma 1**: If an orbiform $\mathcal{O}$ of unit diameter is a subset of $\mathcal{P}(\sigma)$ then it is a subset of $\mathcal{S}(\sigma)$.

Proof: By propositions 3, 4 and 5, $\mathcal{O}$ will not have any points in the interiors of $C_S$, $E_S$ of $A_S$. Therefore $\mathcal{O}$ is a subset of $\mathcal{S}(\sigma)$.

Proof of proposition 6: By proposition 2 any orbiform of unit diameter is congruent to an orbiform $\mathcal{O}$ of unit diameter that is a subset of $\mathcal{P}(\sigma)$. By lemma 1 $\mathcal{O}$ is a subset of $\mathcal{S}(\sigma)$. This proves that $\mathcal{S}(\sigma)$ is a covering for orbiforms of unit diameter and therefore by corollary 1 it is a covering for all shapes of unit diameter.

In the limiting case $\sigma = 0$ the regions $C_S$ and $E_S$ shrink to points and $A_S$ becomes the area removed by Sprague. $\mathcal{S}(0)$ is therefore Sprague's cover. Unfortunately the area of $\mathcal{S}(\sigma)$ is smallest at $\sigma = 0$ so these reductions alone do not provide a new smaller cover. To achieve that, more regions must be removed using reflection arguments similar to those of Hansen. An argument of that type was used previously by Baez, Bagdasaryan and Gibbs to set a new upper bound, but their reduction aimed for simplicity and was not the best possible. Here the objective is to remove a larger area using similar arguments in an attempt to get as close as possible to the optimal covering that can be contained within $\mathcal{P}(\sigma)$.

## Further reductions near corner A

A line drawn from the point $F_3$ to $G$ has length one. Let $\theta$ be the angle that this makes to the edge $ED$ of the hexagon. Consider more generally a line of length one between a point on the side EF to a point $L(s)$ on the side $DC$ which makes an angle s to the side ED with $L(-\theta) = G$. Define also a point $N(t)$ on the edge $DE$ at a distance t from the midpoint $M$ of that edge of the hexagon for $-\tau \leq t \leq \tau$ such that $N(\tau) = E_3$. (see figure 7)

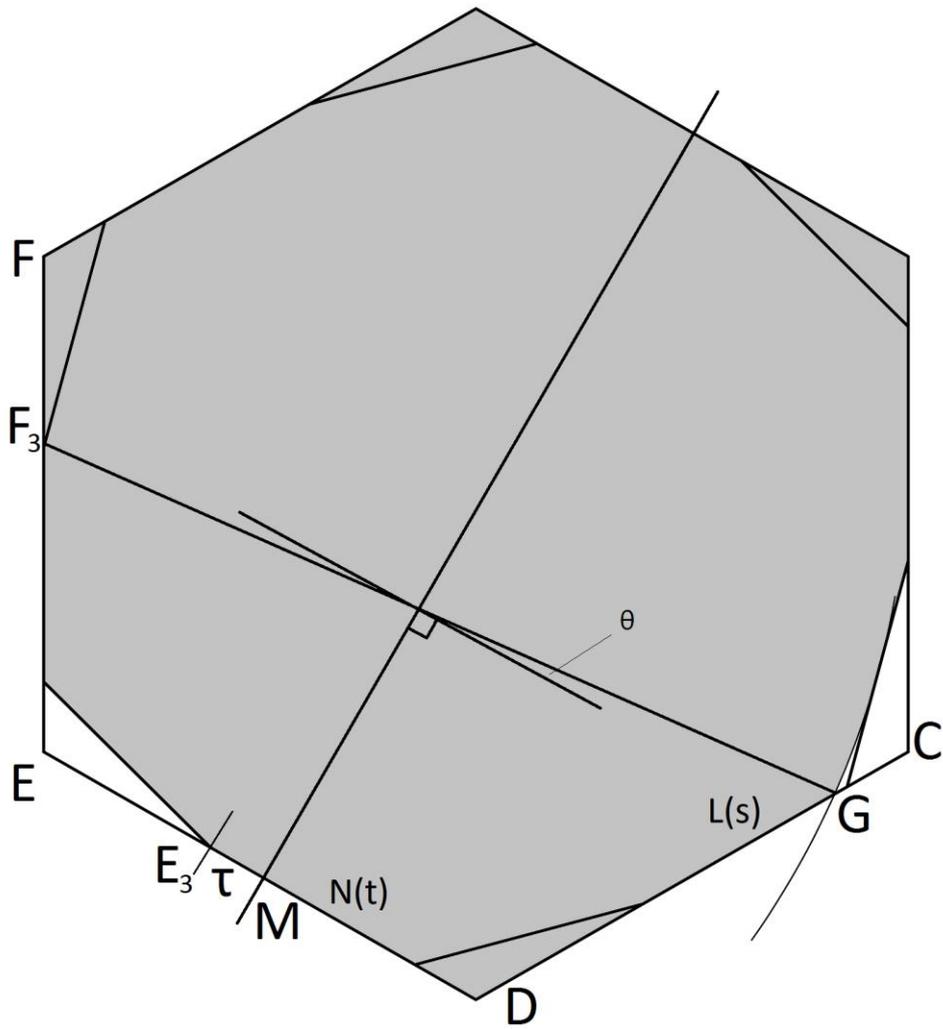

**Figure 7**

Two circular arcs of radius one centred on $L(s)$ and $N(t)$ will meet at a point $X(s,t)$ near the top corner $A$ of the hexagon (figure 8). The set of points $X(s,t)$ for $-\theta \leq s \leq \theta$ and $-\tau \leq t \leq \tau$ is bounded by four arcs. Call this region $\mathcal{R}$. The points $X(-\theta, -\tau) = P$ is on the edge of the hexagon and $X(-\theta, \tau) = X$ as in figure 6. The region $\mathcal{R}$ therefore joins onto two of the arcs that bound the region $A_S$ that has already been removed.

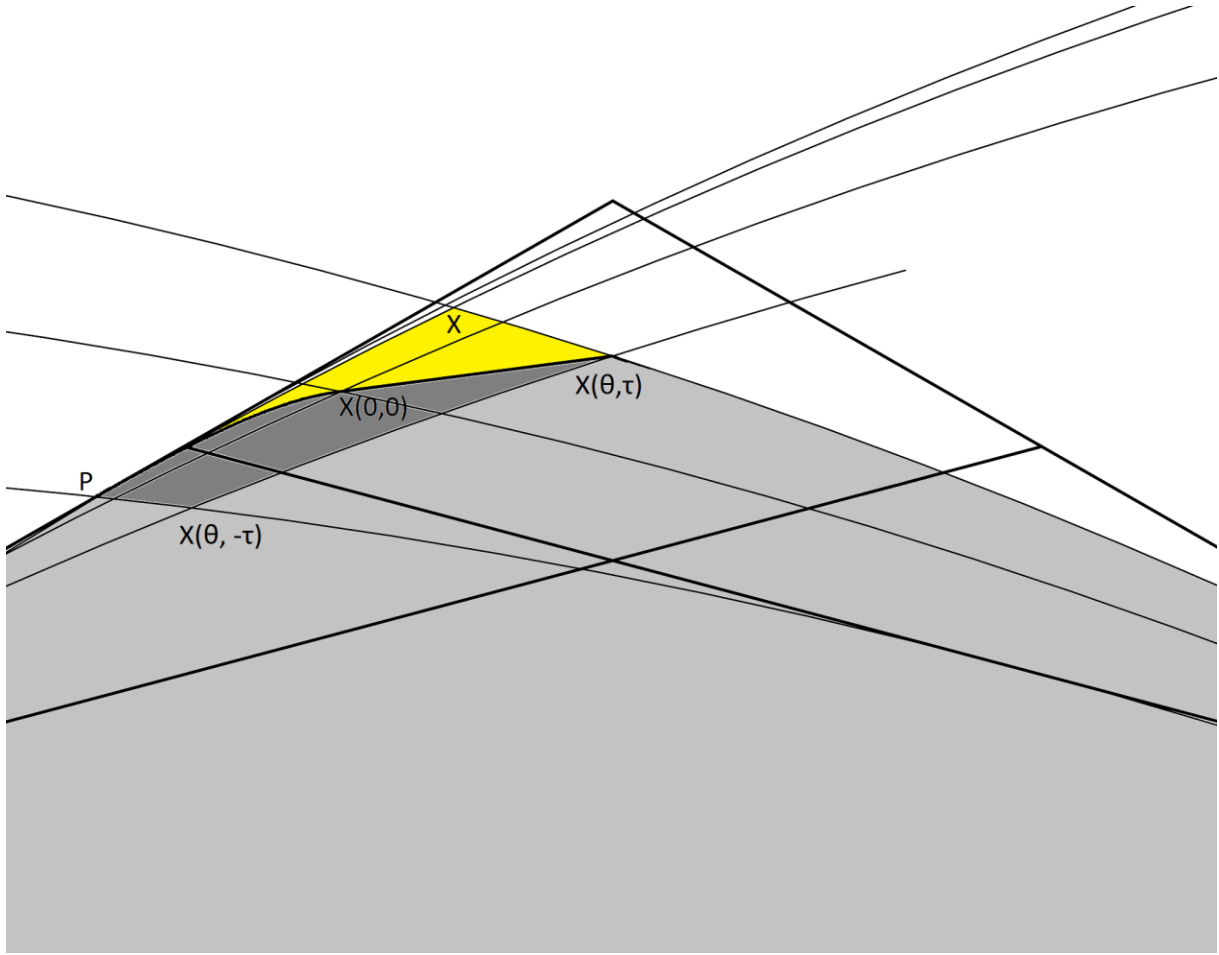

**Figure 8**

In [13] it was shown that part of the region $\mathcal{R}$ with $\tau \geq 0$ can be removed to make a smaller cover. The convex hull of this covering had to be taken according to the original terms of the problem. With a slightly more elaborate argument it is possible to remove a larger portion of $\mathcal{R}$ leaving a smaller covering that is already convex.

**Definition 8**: The **region $A_H$** is defined as follows.

Draw a continuous path from the point $X(\theta, \tau)$ to the point $X(0,0)$, making it the shortest such path that stays within this region $\mathcal{R}$. For large enough slant angles this is a straight line as shown in figure 8, but for smaller slant angles it consists of part of the arc centred on $G$ from $X(\theta, \tau)$ and then a straight line segment tangent to the arc ending at $X(0,0)$. For each point $X(s,t)$ on this path there is an image point $X(-s,-t)$. These points form a curve from $X(0,0)$ to $P = X(-\theta, -\tau)$. The combination of this curve and the path from $X(\theta, \tau)$ to $X(0,0)$ form a curve from $X(\theta, \tau)$ to $P = X(-\theta, -\tau)$. Define a region $A_H \subset \mathcal{R}$ bounded by this curve and the two arcs from $P$ to $X$ and from $X$ to $X(\theta, \tau)$.

**Proposition 7**: Any orbiform of unit diameter is congruent to a subset of $\mathcal{S}(\sigma)$ that does not have any points in the interior of $A_H$.

Take any orbiform $\mathcal{O}$ fitted inside the covering $\mathcal{S}(\sigma)$. To demonstrate that the region $A_H$ can be removed from $\mathcal{S}(\sigma)$ to form a smaller covering it is sufficient to show that if $\mathcal{O}$ enters the interior of

$A_H$ then its reflection about the centre line through M (see figure 7) will also be contained inside $S(\sigma)$ but will not enter the interior of $A_H$.

**Figure 9**

Figure 9 shows the corner regions $A, B, C, D, E, F$ cut at a slant angle $\sigma$ already illustrated in figure 3, but also the regions $A', B', C', D', E', F'$ which are the regions cut at a slant angle $-\sigma$.

**Lemma 2.** If an orbiform $\mathcal{O}$ of unit diameter which is a subset of $S(\sigma)$ contains a point in the interior of region $A_H$ then it has no points in the corner region $D'$.

Proof: All points in the interior of $A_H \subset \mathcal{R}$ are outside the arc of unit radius centered on $N(-\tau)$ and are therefore at a distance greater than one from all points inside the corner region $D'$.

**Lemma 3.** If an orbiform $\mathcal{O}$ of unit diameter which is a subset of $S(\sigma)$ contains a point in the interior of region $A_H$ then it has no points in the corner region $F'$.

Proof: All points in the interior of $A_H \subset \mathcal{R}$ are outside the arc of unit radius centered on $L(\theta)$. Therefore $\mathcal{O}$ must touch the side DC within the corner region $C'$. It therefore cannot enter the region $F'$.

The orbiform $\mathcal{O}$ must therefore be contained within the shaded region shown in figure 9. This tells us that if it is reflected about the centre line through M it will remain within

Proof of proposition 7: Suppose now that the orbiform $\mathcal{O}$ fitted inside $\mathcal{S}(\sigma)$ enters the interior of the region $A_H$. $\mathcal{O}$ must touch the line $ED$ at a unique point $N(t)$ where t is the signed distance from $M$ as defined previously. t must be in the range $-\tau \leq t \leq \tau$ because if $t > \tau$ it would enter the interior of E which is not part of $\mathcal{S}(\sigma)$ and if $t < -\tau$ it would enter $D'$ contrary to lemma 2.

$\mathcal{O}$ must also touch the side $DC$ at a unique point $L(s)$ such that $s \geq -\theta$ because the $\mathcal{O}$ is inside $\mathcal{S}(\sigma)$. Also $s \leq \theta$ since otherwise $\mathcal{O}$ could not enter $\mathcal{R}$.

The point $X(s, t)$ must be in the interior of $A_H$ since otherwise $\mathcal{O}$ could not enter the interior of $A_H$. Therefore the point $X(-s, -t)$ lies inside $\mathcal{R}$ but not inside $A_H$.

Let $\mathcal{O}'$ be the reflection of $\mathcal{O}$ about the centre line through $M$. Since $\mathcal{O}$ does not enter the interior of $D'$ or $F'$ it follows that $\mathcal{O}'$ does not enter the interior of $E$ or $C$. $\mathcal{O}'$ is therefore inside $\mathcal{P}(\sigma)$ and therefore also inside $\mathcal{S}(\sigma)$ by lemma 1.

$\mathcal{O}'$ touches the side $ED$ at the unique point $N(-t)$. It also touches the side $DC$ at a unique point $L(-s')$ where s' < s. Since $X(-s, -t)$ is not in $A_H$ it follows that $\mathcal{O}'$ does not enter $A_H$. It follows that either $\mathcal{O}$ or its reflection $\mathcal{O}'$ is a subset of $\mathcal{S}(\sigma)$ that does not enter the interior of $A_H$. This completes the proof of proposition 7.

By proposition 7 $A_H$ can be removed from $\mathcal{S}(\sigma)$ to provide a smaller covering.

## Further reductions near corner E

The reductions near $A$ are the largest contribution to an improved covering but further reductions near $E$ are also possible. Recall that an area $E_S$ has already been removed around $E_2$ along two arcs $IJ$ and $JH$ (figure 5) A further reduction in this area can be constructed assuming that the slant angle is less than 10 degrees.

**Definition 9**: The **region $E_H$** is defined as follows.

Let $R$ be the intersection of the line segment $F_2 F_3$ with the centre line $F_1 C_1$. An arc of radius 1 centred on $R$ will cut the edge $CB$ at a point $S$. For slant angles less than 10 degrees $S$ will be outside the removed corner area $C$. An arc centred on $S$ intersects the arc $JI$ at a point $T$ (see figure 10.) $T$ is first end point of a line that will be used to remove a further region.

Now consider the edge of the regular dodecahedron that would remove the triangle near $E$ if the slant angle were zero, as done originally by Pal. This line would cut the arc $HJ$ at a point $U$. $U$ is the other end point of a line $TU$, but this may not be a straight line segment. It is the line from $T$ to $U$ that forms a convex hull with the arc centred on $S$ and the point $U$. The region $E_H = TUJ$ bounded by this line and the two arcs $UJ$ and $JT$. This is shown zoomed in figure 11.

At $\sigma = 10°$ the region $E_H$ vanishes and can be defined to be empty for larger slant angles

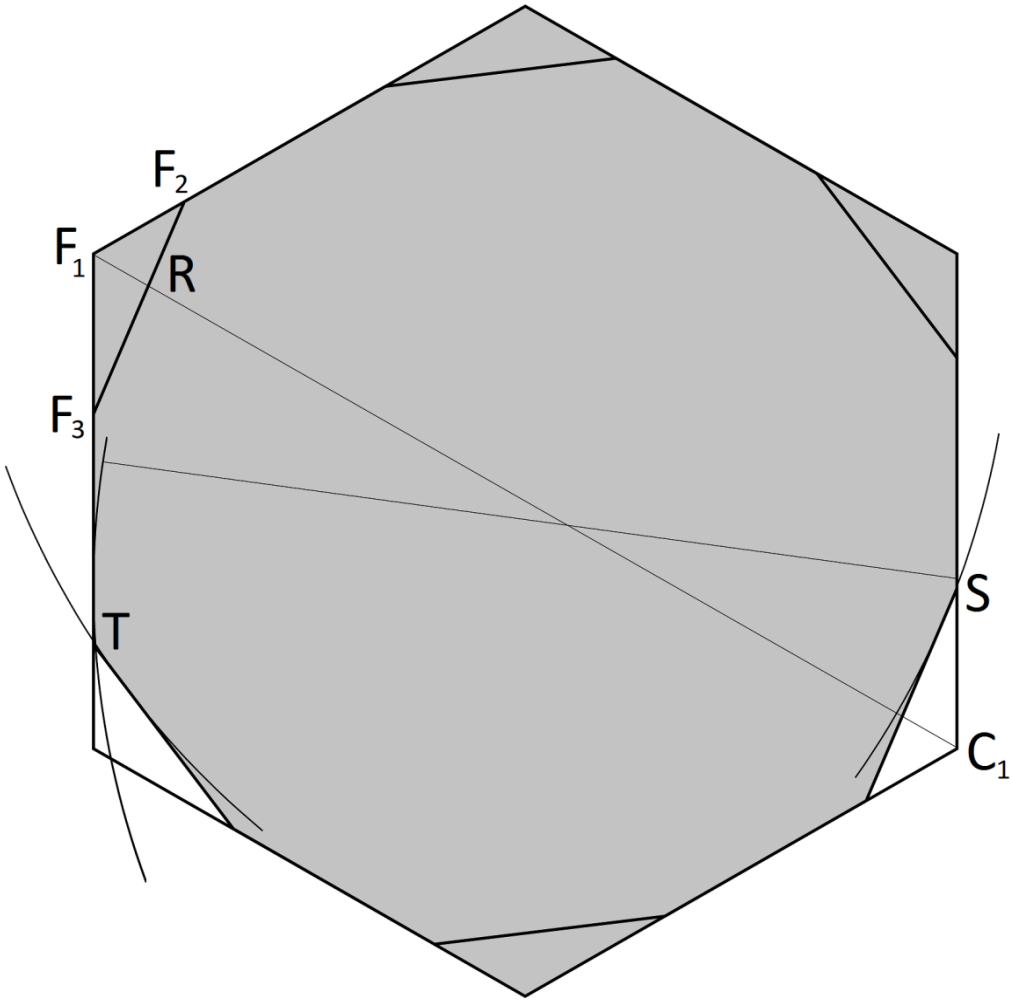

**Figure 10**

**Definition 10**: The **shape $\mathcal{H}(\sigma)$** is defined to be $\mathcal{S}(\sigma)$ with the two areas $A_H$ and $E_H$ removed for $\sigma < 10°$.

**Theorem 1**: $\mathcal{H}(\sigma)$ , $\sigma < 10°$ is a covering for sets of unit diameter.

To prove this consider again an orbiform $\mathcal{O}$ fitted inside $\mathcal{S}(\sigma)$ it is sufficient to show that if $\mathcal{O}$ enters the interior of $A_H$ or $E_H$ then it can be reflected or rotated to a new position in $\mathcal{S}(\sigma)$ so that it does not enter either region.

**Lemma 4.** If an orbiform $\mathcal{O}$ of unit diameter which is a subset of $\mathcal{S}(\sigma)$, $\sigma < 10°$ contains a point in the interior of region $E_H$ then it has no points in the corner region $C'$.

Proof: $\mathcal{O}$ must touch the side $CB$ below the point $S$. This point will be at a distance of more than one from all points on the line segment $F_3R$ but since $\mathcal{O}$ does not enter the interior of $C$ it must cross or touch the line $F_3F_2$ at some point. It follows that this point must be on the segment $RF_2$ and not at $R$, but this is inside the area $F'$. $\mathcal{O}$ therefore has points in the interior of $F'$ and cannot have points in $C'$

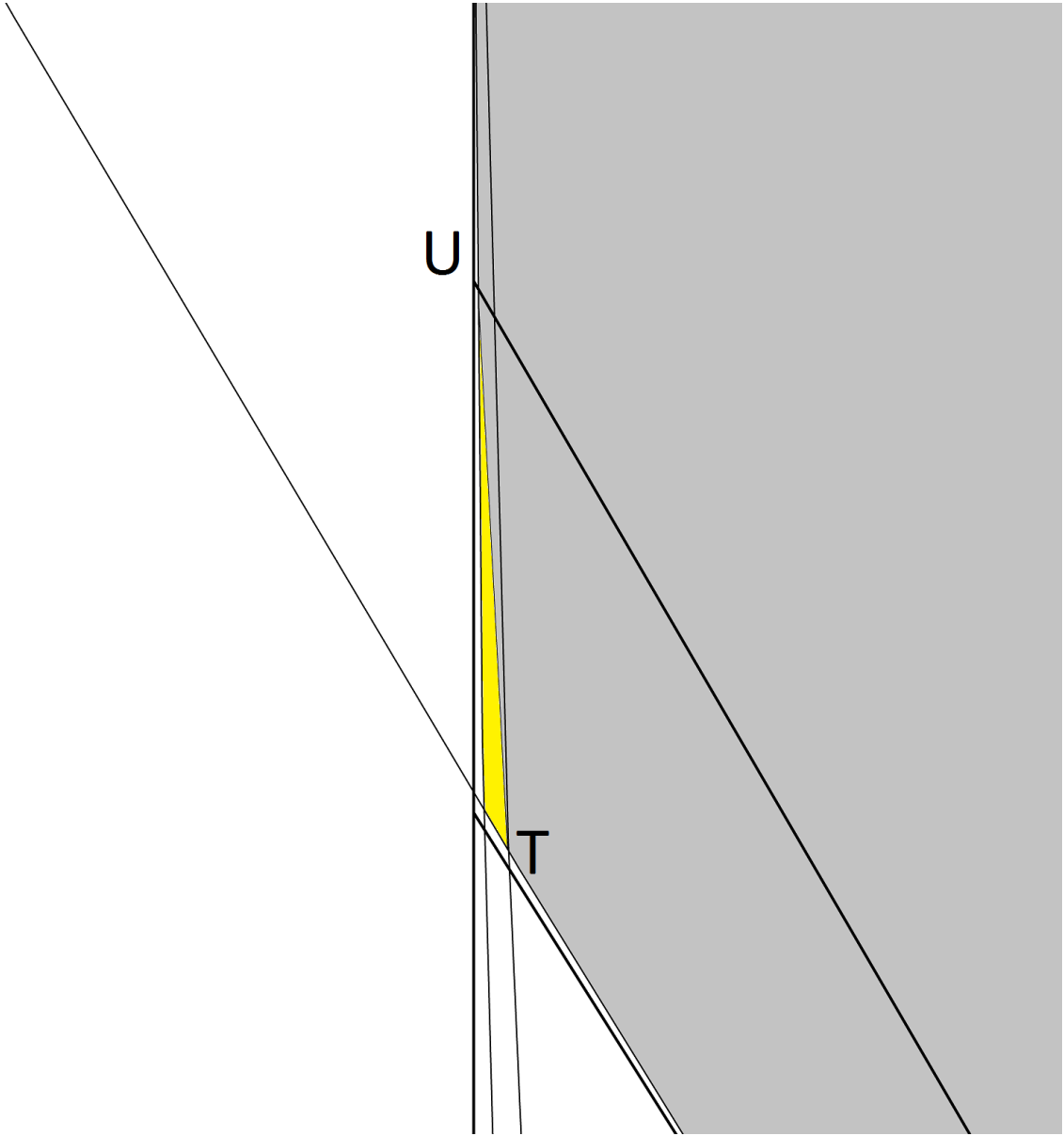

**Figure 11**

**Proposition 8**. If an orbiform $\mathcal{O}$ of unit diameter is a subset of $\mathcal{S}(\sigma)$ for $\sigma < 10°$ then $\mathcal{O}$ cannot enter the interior of both $A_H$ and $E_H$

Proof: This follows from lemma 3 and lemma 4.

**Lemma 5**. If an orbiform $\mathcal{O}$ of unit diameter which is a subset of $\mathcal{S}(\sigma)$, $\sigma < 10°$ contains a point in the interior of region $E_H$ then it has no points in the corner region $B'$.

Proof: $E_H \subset E'$. A point in the interior of $E'$ cannot be in the corner region $B'$

**Proposition 9**. If an orbiform $\mathcal{O}$ of unit diameter which is a subset of $\mathcal{S}(\sigma)$, $\sigma < 10°$ contains a point in the interior of region $A_H$ then it can be reflected into a position where it does not enter the interior of $A_H$ or the interior of $E_H$.

Proof: Use the proof of proposition 7 and observe that the reflected orbifold does not enter the interior of $F'$. Therefore by lemma 4 it does not enter the interior of $E_H$

An orbifold $\mathcal{O}$ of unit diameter that is a subset of $\mathcal{S}(\sigma)$, $\sigma < 10°$ and which has a point in the interior of $E_H$ must by lemma 4 and lemma 5 be within the shaded area shown in figure 12

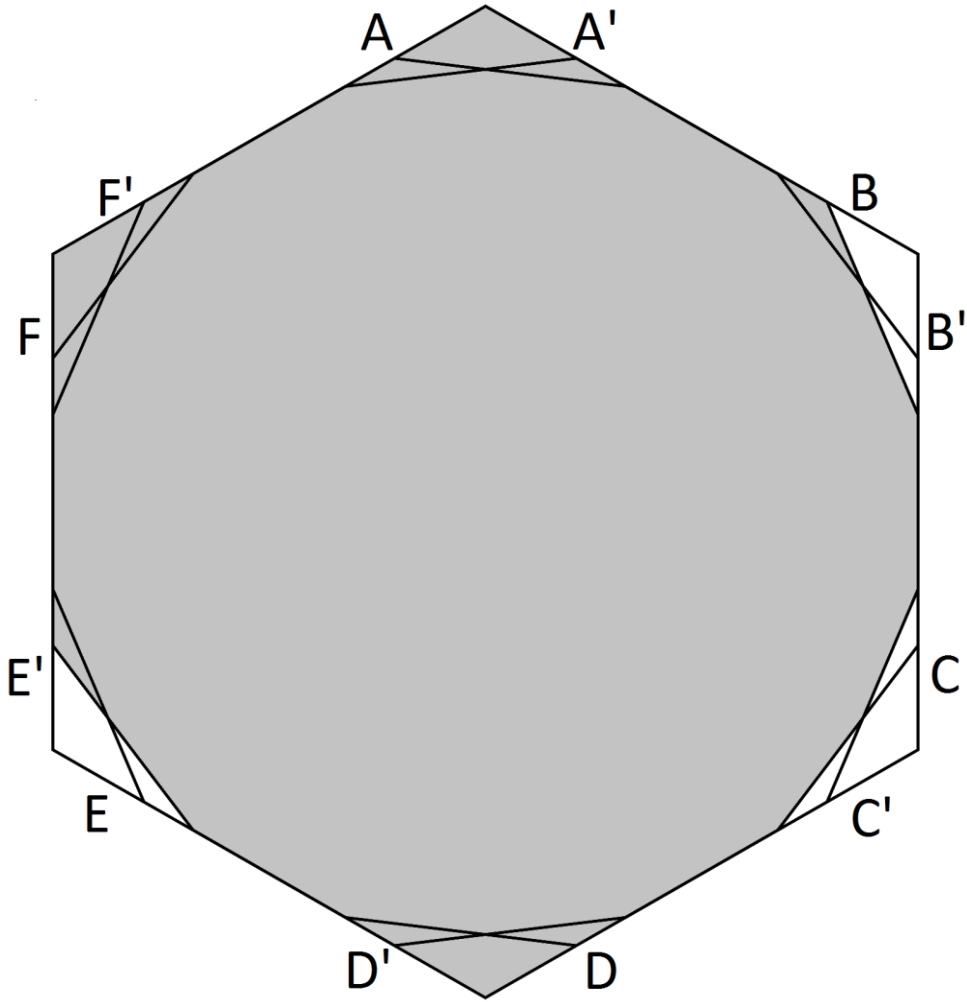

**Figure 12**

The remaining task is to show that in this case $\mathcal{O}$ can always be rotated or reflected so that it does not have any points in the interior of $A_H$ or the interior of $E_H$. To do this the case is divided into three subcases that include all possibilities shown in figure 12 as follows.

Case 1: $\mathcal{O}$ does not enter the interior of $A'$

Case 2: $\mathcal{O}$ does not enter the interiors of $D'$ or $D$.

Case 3: $\mathcal{O}$ does not enter the interiors of $D'$ or $A$.

**Lemma 6** (case 1). If an orbiform $\mathcal{O}$ of unit width is a subset of $\mathcal{S}(\sigma)$, $\sigma < 10°$ and has a point in the interior of $E_H$ but does not enter the interior of $A'$ then its reflection $\mathcal{O}'$ about the centre line through CF will be a subset of $\mathcal{S}(\sigma)$ which does not have a point in the interior of $E_H$ or $A_H$

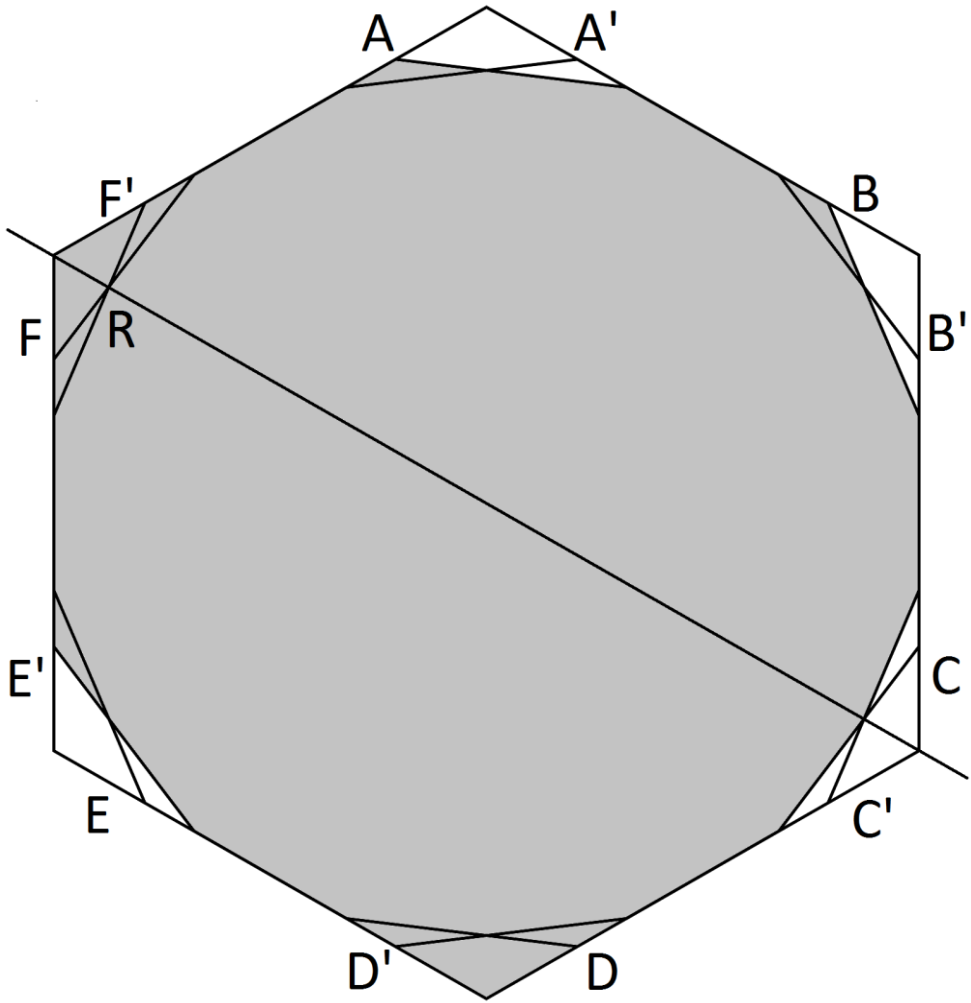

**Figure 13**

Proof: Figure 13 shows the area in grey into which $\mathcal{O}$ can enter according to the premise of case 1. If $\mathcal{O}'$ is the image of $\mathcal{O}$ under reflection in the centre line through $R$ then $\mathcal{O}'$ is also a subset of $\mathcal{S}(\sigma)$.

Since $\mathcal{O}$ enters the interior of $E_H$ it cannot have any points on the line segment $F_3R$ (see proof of lemma 4) but it must enter or touch the corner region F' somewhere. This means that $\mathcal{O}'$ must have a point on the line segment $F_3R$ and therefore it cannot enter the interior of $E_H$. It also cannot enter the interior of $A_H$ by lemma 3.

**Lemma 7** (case 2). If an orbiform $\mathcal{O}$ of unit width is a subset of $\mathcal{S}(\sigma)$, $\sigma < 10°$ and has a point in the interior of $E_H$ but does not enter the interiors of $D$ or $D'$ then its reflection $\mathcal{O}'$ about the centre line through the midpoints of $AF$ and $CD$ will be a subset of $S(\sigma)$ which does not have a point in the interior of $E_H$ or $A_H$

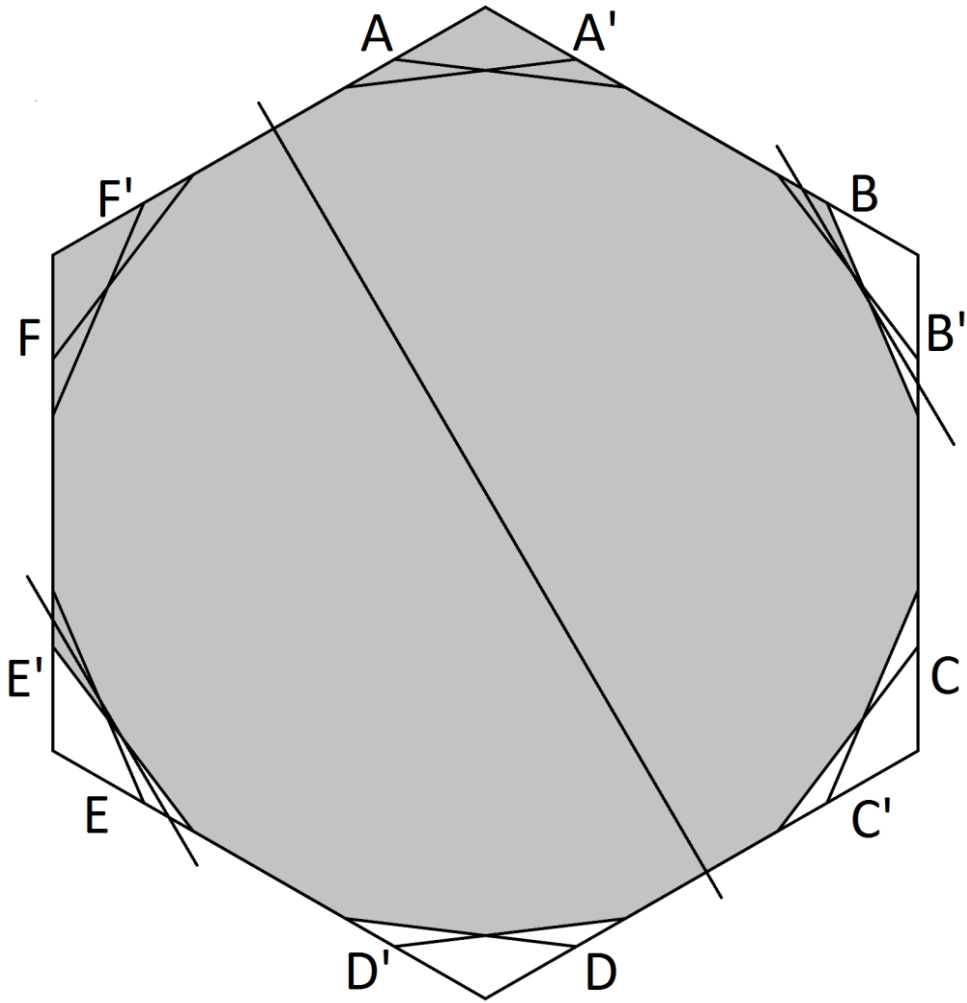

**Figure 14**

Proof: Figure 14 shows the area in grey into which $\mathcal{O}$ can enter according to the premise of case 2. If $\mathcal{O}'$ is the image of $\mathcal{O}$ under reflection in the midpoints of $AF$ and $CD$ then $\mathcal{O}'$ is also a subset of $S(\sigma)$.

Draw also two lines parallel to the centre line at a distance of one half to either side. By the definition of $E_H$, $\mathcal{O}$ must fall outside the line that passes near corner E. This means that $\mathcal{O}'$ falls outside the other line where it passes near B. Since the two parallel lines are separated by unit distance $\mathcal{O}'$ cannot enter the interior of region $E_H$. It also cannot enter the interior of $A_H$ by lemma 3.

**Lemma 7** (case 3). If an orbiform $\mathcal{O}$ of unit width is a subset of $S(\sigma)$, $\sigma < 10°$ and has a point in the interior of $E_H$ but does not enter the interiors of $A$ or $D'$ then its image $\mathcal{O}'$ under a rotation of 120° about the centre point will be a subset of $\mathcal{S}(\sigma)$ which does not have a point in the interior of $E_H$ or $A_H$

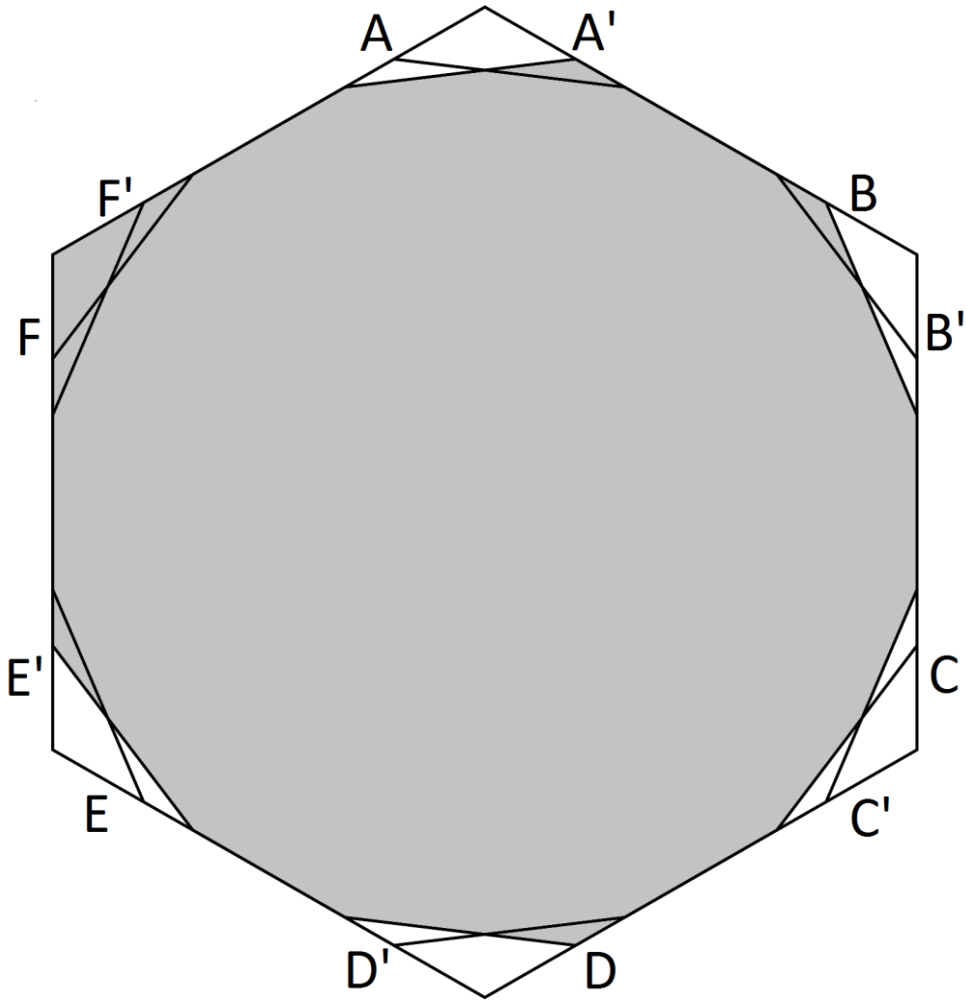

**Figure 15**

Proof: Figure 15 shows the area in grey into which $\mathcal{O}$ can enter according to the premise of case 3. A rotation of 120° clockwise maps $A$ to $E$ and $E$ to $C$. Therefore the image $\mathcal{O}'$ of $\mathcal{O}$ under the rotation is also a subset of $\mathcal{S}(\sigma)$.

Furthermore $\mathcal{O}'$ cannot enter the interior of $A_H$ by lemma 2 and it cannot enter the interior of $E_H$ by lemma 4

**Proposition 10.** If an orbiform $\mathcal{O}$ of unit width is a subset of $\mathcal{S}(\sigma)$, $\sigma < 10°$ and has a point in the interior of $E_H$ then there is an orbiform $\mathcal{O}'$ congruent to $\mathcal{O}$ which is a subset of $\mathcal{S}(\sigma)$ and which does not enter the interior of $E_H$ or $A_H$

Proof: At least one of the cases 1,2 or 3 must apply therefore by lemmas 5,6 and 7 there must be an image $\mathcal{O}'$ of $\mathcal{O}$ under a reflection or rotation that is a subset of $\mathcal{S}(\sigma)$ and which does not enter the interior of $E_H$ or $A_H$

Proof of theorem 1: By proposition 6 $\mathcal{S}(\sigma)$ has a subset $\mathcal{O}$ congruent to any orbiform of unit diameter. By propositions 8, 9 and 10 any there is an orbiform $\mathcal{O}'$ congruent to $\mathcal{O}$ that does not enter the interior of $E_H$ or $A_H$. Therefore by definition 10 $\mathcal{H}(\sigma)$, $\sigma < 10°$ is a covering for sets of unit diameter.

## Upper bounds for Lebesgue's covering problem

$\mathcal{H}(\sigma)$ is convex by construction and for $\sigma < 10°$ is a covering for sets of unit diameter by theorem 1. The area of $\mathcal{H}(\sigma)$ for any value of $\sigma < 10°$ therefore sets an upper bound for the answer to Lebesgue's covering problem. To find the best bound the minimum area should be computed. Figure 16 shows a plot of the area of $\mathcal{H}(\sigma)$ for $\sigma < 4°$.

The computed minimum area is less than 0.8440935944 at $\sigma < 1.5494°$.

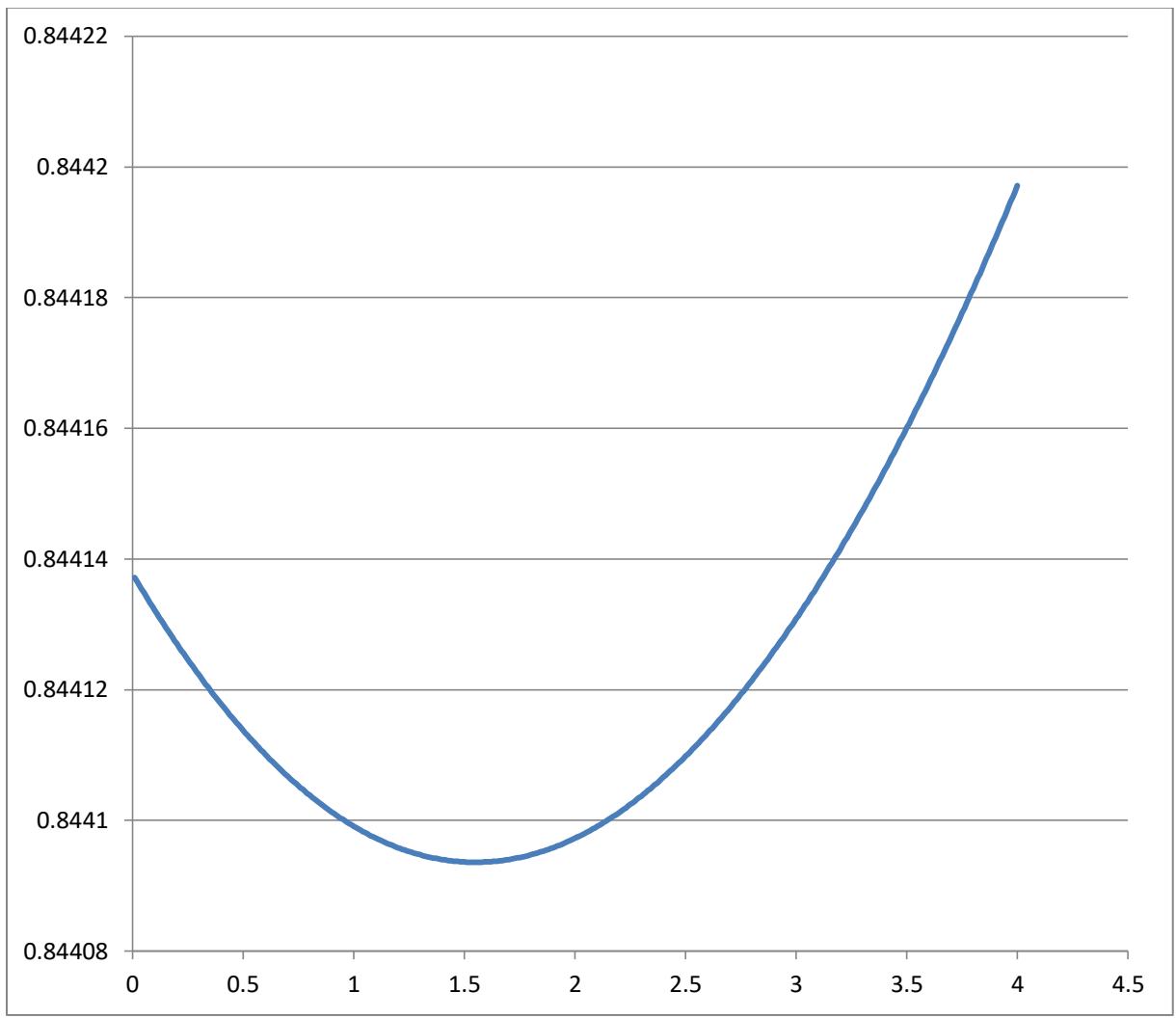

**Figure 16**

## Future prospects

What can be said about the solution to Lebesgue's covering problem in the light of this research? Setting lower bounds has proved harder than setting upper bounds. However, computational studies provide some evidence that the optimal covering is a subset of $\mathcal{P}(\sigma)$ for some angle $\sigma$ [12]. If this is assumed as a hypothesis then conditional lower bounds can be found by looking for orbiforms that are congruent to a unique subset of $\mathcal{P}(\sigma)$ for each $\sigma$, and taking the area of their convex hull. This can then be minimised to give a conditional lower bound. A necessary but not sufficient condition for an orbiform to have a unique fit in $\mathcal{P}(\sigma)$ is that it has an axis of reflection symmetry. A computational check can be performed on any given shape to see if its fit is unique. In this way it has been found that only the regions near the existing reductions $E_H$ or $A_H$ have further scope for reducing the covering area further. Indications are that such further reductions are indeed possible at both locations.

## Acknowledgments

I am grateful to John Baez, Karine Bagdasaryan and Greg Egan for useful discussions and input to this problem.

# References

[1] Henri Lebesgue (1914), Letter to Julius Pal, described in [4]

[2] Heinrich W.E. Jung, "Über die kleinste Kugel die eine räumliche Figur einschliesst", *Reine Agnew Math* **123** (1901) 241-257

[3] Heinrich W.E. Jung, "Über den kleinsten Kreis der eine ebene Figur einschliesst", *Reine Agnew Math* **137** (1910) 310-313

[4] Julius (Gyula) Pál, "Über ein elementares Variationsproblem", *Danske Mat.-Fys. Meddelelser* III. 2. (1920)

[5] Roland P. Sprague, "Über ein elementares Variationsproblem". *Matematiska Tidsskrift* Ser. B: 96–99, (1936)

[6] Herbert Meschkowski, Ungelöste und unlösbare Probleme der Geometrie, Braunschweig 1960, (translated as Unsolved and unsolvable problems in geometry, Oliver and Boyd 1966, translated by Jane A.C. Burlak)

[7] H.C. Hansen, 'A small universal cover of figures of unit diameter', *Geom. Dedicata* **4**, 165-172. (1975)

[8] H.C. Hansen, "Towards the minimal universal cover", *Normat* **29** 115-119, 148, (1981)

[9] H.C. Hansen, "Small Universal covers for sets of unit diameter", *Geometriae Dedicata* **42** (1992) 205-213, (1992)

[10] György Elekes, "Generalized Breadths, Circular Cantor Sets, and the Least Area UCC", *Discrete & Computational Geometry*, Volume 12, Issue 1, 439-449, (1994)

[11] Peter Brass, Mehrbod Sharifi, "A Lower Bound for Lebesgues' Universal Covering Problem", *Int. J. Comput. Geom. Appl*. 15, 537, (2005)

[12] Philip E. Gibbs, "A New Slant on Lebesgue's Universal Covering Problem", arXiv:1401.8217, (2014)

[13] John C. Baez, Karine Bagdasaryan, Philip E. Gibbs, "The Lebesgue Universal Covering Problem", Journal of Computational Geometry **6**(1), 288–299, (2015)

[14] Paul J. Kelly, Max L. Weiss, "Geometry and Convexity: A Study in Mathematical Methods." Wiley. pp. Section 6.4, (1979).

[15] M.D. Kovalev, "A Minimal Lebesgue Covering Exists", Mathematical notes of the Academy of Sciences of the USSR, Volume 40, Issue 3 (1986) pp 736–739, (1985)

[16] S. Vrecica, "A note on curves/sets of constant width", *Publications de L'Institut Mathematique* **29**, 289–291, (1981),